\documentclass[12pt]{article}
\usepackage{graphicx,subfigure}
\usepackage{color,latexsym,amsfonts,amssymb}
\usepackage{bm}     
\usepackage{amsthm} 
\usepackage{amsmath}
\usepackage{booktabs} 
\usepackage{boxedminipage}  
\usepackage{natbib} 
\usepackage{algorithm,algpseudocode}  
\usepackage{subfigure}

\newtheorem{theorem}{Theorem}

\newtheorem{lemma}{Lemma}

\newtheorem{assumption}{Assumption}

\newtheorem{example}{Example}

\newcommand{\cU}{\mbox{$\cal U$}}

\newcommand{\cS}{\mbox{$\cal S$}}

\newcommand{\ep}{\epsilon}
\newcommand{\be}{\begin{equation}}
\newcommand{\ee}{\end{equation}}
\newcommand{\ba}{\begin{eqnarray}}
\newcommand{\ea}{\end{eqnarray}}
\newcommand{\bas}{\begin{eqnarray*}}
\newcommand{\eas}{\end{eqnarray*}}
\newcommand{\te}{\hfill $\Box$}

\title{On the existence of optimal stationary policies for average Markov decision processes with countable states}
\author{Li Xia, { Xianping Guo}, { Xi-Ren Cao\thanks{L. Xia
is with the Business School, Sun Yat-Sen University, Guangzhou
510275, China (email: xiali5@sysu.edu.cn). X. Guo is with the School
of Mathematics, Sun Yat-Sen University, Guangzhou 510275, China
(email: mcsgxp@mail.sysu.edu.cn). X.-R. Cao is with the Shanghai
Jiao Tong University and also with the Hong Kong University of
Science and Technology (email: eecao@ust.hk. \emph{Corresponding
Author}).}}}

\date{}

\begin{document}

\maketitle

\begin{abstract}
For a Markov decision process with countably infinite states, the
optimal value may not be achievable in the set of stationary
policies. In this paper, we study the existence conditions of an
optimal stationary policy in a countable-state Markov decision
process under the long-run average criterion. With a properly
defined metric on the policy space of ergodic MDPs, the existence of
an optimal stationary policy can be guaranteed by the compactness of
the space and the continuity of the long-run average cost with
respect to the metric. We further extend this condition by some
assumptions which can be easily verified in control problems of
specific systems, such as queueing systems. Our results make a
complementary contribution to the literature in the sense that our
method is capable to handle the cost function unbounded from both
below and above, only at the condition of continuity and ergodicity.
Several examples are provided to illustrate the application of our
main results.
\end{abstract}
\textbf{Keywords}: Markov decision process, countable states,
optimal stationary policy, metric space

\section{Introduction}\label{section_intro}
For finite Markov decision processes (MDPs), the optimality of
various types of policies are well studied. For example, it is well
known that the optimal value of finite MDPs with discounted or
average criteria can be achieved by \emph{Markovian} and
\emph{deterministic} policies, thus \emph{history-dependent} and
\emph{randomized} policies are not needed to consider. More details
can be referred to books on MDPs \citep{Bertsekas12,Puterman94}.

Countable-state MDPs are a type of widely existing models and are
particularly useful for many problems, such as queueing systems,
inventory management, etc. When the state space of MDPs is changed
from finite to infinite (countable), the relevant analysis becomes
more complicated and the algorithms need sophisticated discussion
\citep{Golubin2003,Meyn1997}. Compared with the complete theoretical
results for finite MDPs, there is no comprehensive theory for
infinite MDPs with countable states and the long-run average
criterion. The existence of an optimal stationary policy for
countable-state MDPs needs specific discussion, and attracts
research attention in recent decades. Although we can restrict our
attention to stationary policies in finite MDPs, this is no longer
true when the state space is countable. In general, the optimal
value of a countable-state MDP may not be achievable by stationary
policies, even not by history-dependent policies. Interesting
counterexamples can be found in the excellent books on MDPs (see
Examples 5.6.1\&5.6.5\&5.6.6 of \cite{Bertsekas12}, Examples
8.10.1\&8.10.2 of \cite{Puterman94}, and Subsection~7.1 of
\cite{Sennott1999}).

Since a stationary policy is not necessarily optimal for countably
infinite MDPs, there are literature works on the specific existence
conditions of optimal stationary policies. Sennott studies the
existence conditions for average cost optimality of stationary
policies for discrete-time MDPs when state space is countable and
action space is finite \citep{Sennott1986,Sennott1989}. In Sennott's
studies, a distinguished state is introduced and the vanishing
discount optimality approach is adopted to study the optimality
inequality. \cite{Borkar1989} also studies the condition of optimal
stationary policies for discrete-time average cost MDPs with
countable states, but from the characterization through the dynamic
programming equations. For constrained MDPs with countable states
and long-run average cost, \cite{Borkar1994} further establishes the
existence of stationary randomized policies for the general case of
nonnegative cost functions (or unbounded from below), which uses the
method of occupation measures. \cite{Lasserre1988} studies the
stationary policies of denumerable state MDPs for not only the
average cost optimality, but also the Blackwell optimality.
\cite{Meyn1999} studies the similar problem based on the
stabilization of controlled Markov chains with algorithmic analysis.
\cite{Cao2015} study the existence condition of optimal stationary
policies for a class of queueing systems, also from the analysis of
system stability. \cite{Cavazos1991,Cavazos1992} give a fairly
complete summary and comparison of different results on existence
conditions for discrete-time average cost MDPs with countable state
space and finite action sets.

For more general cases rather than countable state space,
\cite{Hernandez1991} studies the existence condition on average cost
optimal stationary policies in a class of discrete-time Markov
control processes with Borel spaces and unbounded costs, where the
action space is assumed setwise continuity instead of a compact set.
\cite{Feinberg2007} present sufficient conditions for the existence
of an optimal stationary policy of MDPs with the average cost
optimality inequalities, where the state and action space are Borel
subsets of Polish spaces. The derived result is also applied to a
cash balance problem with an inventory model. For continuous-time
MDPs with infinite state in Polish spaces, \cite{Guo2006} study the
existence of optimal deterministic stationary policies by using the
Dynkin formula and two optimality inequalities for the average cost
criterion. Some other systematic discussion on this issue can also
be found in the excellent books on MDPs, see
\cite{Bertsekas12,Hernandez1996,Puterman94,Sennott1999} for
discrete-time MDPs and \cite{Bertsekas12,Guo2009} for
continuous-time MDPs.

In summary, most of the existing results are about the sufficient
conditions, which usually require constructing a set of functions
satisfying several sophisticated assumptions. Although these
conditions are quite general, they may be not easy to verify and may
encounter difficulty of function construction during the application
to practical problems. In this paper, we study the optimality
condition of stationary policies for average cost MDPs with
countable states and finite actions available at each state. By
defining a proper metric in the policy space, we study the
continuity of the system's average cost and the compactness of the
policy space, and we show that such continuity and compactness can
induce the existence of an optimal stationary policy. We further
extend the continuity requirement by assuming some reasonable
conditions on transition rates and uniform convergence of
un-normalized probabilities in MDPs. Compared with the existing
literature work, our result holds at a weak condition of requiring
continuity and ergodicity, and it can handle the cost function
unbounded from both below and above. While some general results in
the literature require the cost function unbounded only from below
(e.g., see \citep{Borkar1994}) or $\omega$-geometric ergodicity
(e.g., see \citep{Hernandez1999}), which partly demonstrates the
advantages of our method. Moreover, our result may be easier to
verify for some MDPs, especially for queueing systems. The main
results of the paper are illustrated by several examples, for one of
which the cost function is unbounded from above and from below, as
discussed in Remark~2 at the end of Section~\ref{section_queue}.

The remainder of the paper is organized as follows. In
Section~\ref{section_result}, we derive the existence condition by
studying the continuity of the average cost in a defined compact
metric space of policies. In Section~\ref{section_queue}, an example
of scheduling problem in queueing systems is provided to demonstrate
the validation process of our existence condition of an optimal
stationary policy. In Section~\ref{section_extension}, we further
extend the existence condition to several reasonable assumptions
which may be easy to satisfy in practical problems. Finally, we
conclude the paper in Section~\ref{section_conclusion}.

\section{The Basic Idea}\label{section_result}
In an MDP, the state space is denoted as $\mathcal S$, which is
assumed to be countably infinite. Without loss of generality, we
denote it as $\mathcal S =\{0,1,\dots\}$. Associated with every
state $i \in \mathcal S$, there is a finite action set $\mathcal
A(i)$. At state $i \in \mathcal S$, if action $a \in \mathcal A(i)$
is adopted, an instant cost $f(i,a)$ will incur. Meanwhile, the
system will transit to state $j \in \mathcal S$ with transition
probability $p^a(i,j)$ for discrete-time MDPs and with transition
rate $q^a(i,j)$ for continuous-time MDPs, respectively. Let $u$
denote a (deterministic) stationary policy which is a mapping on
$\mathcal S$ such that $u(i) \in \mathcal A(i)$ for all $i \in
\mathcal S$. Let $\cU$ denote the stationary policy space and $\cU:=
\times_{i \in \mathcal S}\mathcal A(i) := \mathcal A(0) \times
\mathcal A(1) \times \dots $, with ``$\times$" being the Cartesian
product. Let $X (t)$ be the system state at
time $t$. Under suitable conditions, the long-run average
performance measure for MDPs, which does not depend on any initial
state $x\in \cS$, but depends on $u \in \cU$, is defined as $\eta
(u)$:
\begin{equation}
\eta (u) := \lim\limits_{T\rightarrow \infty}\frac{1}{T} \mathbb
E\left\{\sum_{t=0}^{T-1} f(X(t),u(X(t))) \Big | X(0)=x \right\},
\end{equation}
or
\begin{equation}
\eta (u) := \lim\limits_{T\rightarrow \infty}\frac{1}{T} \mathbb
E\left\{\int_{t=0}^{T} f(X(t),u(X(t)))dt \Big | X(0)=x \right\},
\end{equation}
for discrete-time and continuous-time ergodic MDPs, respectively,
where the expectation operator  $\mathbb E$ depends on $u \in \cU$.
However, such dependence is omitted below for notation
simplicity. The goal of optimization is to find a policy $u^*$ such
that
\begin{equation}
\eta (u^* ) = \inf_{u \in \cal U} [\eta (u) ],~~~~~(\mbox{or }~
 \eta (u^* ) = \sup_{u \in \cal U} [\eta (u) ] ).
\end{equation}
Assume that $\eta (u)$ is bounded in $u \in \cU$, so $\inf_{u \in
\cal U} [\eta (u) ]$  is finite. We aim to find conditions under
which such an optimal stationary policy $u^*$ exists.

\begin{theorem}  \label{thmopex}
Suppose $\cU$ is a compact metric space and the function $\eta (u)$
is continuous in $\cU$ with the metric, then an optimal policy $u^*$
exists.
\end{theorem}
{\it Proof:} Let $\eta^* := \inf_{u \in \cal U} [\eta (u) ]$. By definition,
there exists a sequence of policies, denoted as $u_0$, $u_1$, $\dots$,
such that \be     \label{etatos} \lim_{n \to \infty} \eta (u_n ) =
\eta^* . \ee
Because $\cU$ is  compact, there is a
subsequence of $\{ u_n , n=0,1, \dots \}$ that converges to a limit
(accumulation) point. Denote this subsequence as $\{ u_{n_k}, k=0,1,
\dots \}$ and the limit point as $u^* \in \cU$. Then
 \[
 \lim_{k \to \infty } u_{n_k} = u^* \in \cU .
 \]
By continuity of $\eta (u)$, we have
 \[
 \lim_{k \to \infty } \eta (u_{n_k} ) = \eta (u^* ) .
 \]
By (\ref{etatos}), we obtain
\[
\eta (u^*) = \eta ^* =  \inf_{u \in \cal U} [\eta (u) ] ;
\]
i.e., $u^* \in \cU$ is an optimal policy.   \te

Theorem~\ref{thmopex} requires a compact metric space defined for
$\mathcal U$. Below, we introduce such a metric in the policy space.
Note that a policy can be denoted as
 \[
 u = ( u(0) , u(1) , \dots ).
 \]
Choosing a real number $0<r<0.5$, (e.g., $r=0.1$), we define the
distance between two policies $u_1 = ( u_1 (0) , u_1 (1) , \dots )$
and $u_2 = ( u_2 (0) , u_2 (1) , \dots ) $ as
\be   \label{defdis}
  d(u_1 , u_2 ) := \sum_{i=0}^\infty ||u_1 (i) - u_2 (i) || r^i ,
\ee
in which
\[
  ||u_1 (i) - u_2 (i) || := \left \{
  \begin{array}{ll}
  1  &  ~ if ~ u_1 (i) \neq u_2 (i) , \\
  0 &   ~ if ~ u_1 (i) = u_2 (i) .
 \end{array}
 \right .
\]
It is easy to verify that
\[
d(u,u)=0,~~ d(u_1, u_2 ) = d(u_2 , u_1 ),
\]
and for any three policies $u_1$, $u_2$, and $u_3$,
 the following triangle inequality holds
 \[
 d(u_1, u_3 ) \leq d(u_1 , u_2 ) + d(u_2 , u_3 ).
 \]
Thus, $d(u_1 , u_2 )$, $u_1, u_2 \in \cU$, indeed defines a metric
on $\cU$.

Suppose for two policies $u_1$ and $u_2$, $u_1 (i) = u_2 (i)$ for
all $i=0,1, \dots, k$. Then
\begin{eqnarray}\label{drskm}
d(u_1 , u_2 ) &=& \sum_{i=k+1}^\infty ||u_1 (i) - u_2 (i) || r^i
\nonumber\\
&\leq& \sum_{i=k+1}^\infty r^i = r^{k+1} \sum_{i=0}^\infty r^i =
\frac {r^{k+1}}{1-r} < r^{k},
\end{eqnarray}
where the last inequality holds because we choose $r<0.5$, so $\frac
{r}{1-r} <1$. By (\ref{drskm}), we have
\begin{lemma}  \label{lemduout}
$d(u_1, u_2) < r^k$ if and only if $u_1 (i) = u_2 (i)$ for all $i
\leq k$.
\end{lemma}
{\it Proof:} The ``If" part follows directly from (\ref{drskm}). Now
we prove the ``Only if" part using contradiction. Assume that there
is an integer $n$ such that $u_1 (n) \neq u_2 (n)$ and $n \leq k$.
By (\ref{defdis}), we have $d(u_1, u_2) \geq r^n >r^k$, which is in
contradiction with the condition $d(u_1, u_2) < r^k$. Thus, the
assumption is not true and the ``Only if" part is proved. \te

The metric defined by the distance function $d(u_1, u_2 )$ induces a
topology on $\cU$. First, we define an open ball around a point $u
\in \cU$ as
\begin{equation} \label{eq_disc}
 O_\ep (u) := \{ all~v \in \cU: d(u,v) < \ep \} ,~~~\ep >0 .
\end{equation}
We have $u \in O_\ep (u)$ for any $\ep >0$. A set $N(u)$ is called a
neighborhood of a point $u \in \cU$, if there is an open ball $O_\ep
(u)$ for some $\ep >0$ such that $O_\ep (u) \subseteq N(u)$.

By Lemma~\ref{lemduout}, we have the following fact:  $u'(i)=u(i)$
for all $i \leq k$ if and only if $u' \in O_{r^k}(u)$.

\noindent\textbf{Remark~1.} Lemma~\ref{lemduout} reveals the
advantage of the metric (\ref{defdis}): It shows that all the
policies in a small neighborhood $O_{r^k}(u)$ of policy $u$ take the
same actions in the first $k$ states. This property is very useful
in proving the continuity of $\eta (u)$ in many optimization
problems, in which the steady-state probability of state $i$, $\pi
(i)$, goes to zero when $i$ goes to infinity; in other words,
states $i>k$ are less important. \te

In a metric space $\cU$, a limit point can be defined by the metric, i.e.,
$\lim_{n \to \infty} u_n =u$ for some sequence $\{u_n\} \subseteq \mathcal U$,
if and only if $\lim_{n \to \infty} d(u_n , u) =0$. In
this sense, a continuous function is defined in the same way as a
continuous function defined in a real space.

Since $\mathcal A(i)$ is finite and $\mathcal S$ is countable, it is
well known that with the metric (\ref{defdis}) the policy space
$\cU=\times_{i \in \mathcal S}\mathcal A(i)$ is compact. In fact, every
point $u \in \cU$ is an accumulation (limit) point, and every policy is
in $\cU$. In order to apply Theorem~\ref{thmopex}, we have to prove the
continuity of $\eta (u)$ in $\cU$ for the specific problems. Below,
we use some examples to illustrate the applicability of
Theorem~\ref{thmopex} in MDPs.
%

\begin{example} \label{exmp} (A modification of Example 8.10.2 in Puterman's
book \citep{Puterman94})
 Consider an MDP with $\cS = \{1,2, \dots\}$. At each state $i \in \cS$,
 there are two actions $1$ and $0$. If action $1$ is taken,
 then the state transits from $i$ to $i+1$ with probability $1$ and the cost is $f(i,1)=0$;
 if action $0$ is taken, then the state stays at
 $i$ with probability $1$ and the cost is $f(i,0)=\frac {1}{i}$. The
 Markov chain (under any given policy) is denoted as $X(t)$, $t=0,1, \dots$. A stationary policy
 is denoted as a mapping $u: \cS \rightarrow \{0,1\} $.

 The performance measure for policy $u = ( u(1), u(2),
 \dots )$ with initial state $i$ is the long-run average
  \be    \label{defperf}
  \eta (u,i) = \lim_{T \to \infty} \frac {1}{T} \mathbb E \left
  \{\sum_{t=0}^{T-1}
  f(X(t),u(X(t))) \Big | X (0) = i \right \} .
  \ee
 Note that the performance may depend on the initial state, i.e., it
 is a function of both the initial states and policies. To prove the
 existence of an optimal policy, we need to fix the initial state. In
 (\ref{defperf}), we choose $X(0)=1$. We wish to find a policy $u^*$
 such that
 \[
 \eta (u^*, 1) = \inf_{u \in \cal U} \{ \eta (u,1) \} .
 \]
 We need to prove that such an optimal stationary policy exists.

 Now, we prove that $\eta (u,1)$ is continuous in $u \in \cU$ with metric (\ref{defdis}).
 Given a policy $u_0$, for any small positive $\ep$, we find the
 maximum $k$ satisfying $r^k > \ep$. By Lemma \ref{lemduout}, if we choose
 a policy $u$ satisfying $d(u,u_0 )<\ep$, then all the actions of such policies
 $u$ and $u_0$ at states $i \leq k$ are the same. By the structure of $\eta (u,i)$
 defined in (\ref{defperf}), we can conclude that
 \[
 |\eta (u,1) - \eta (u_0,1)| < \frac {1}{k} .
 \]
 More precisely, since $u(i)=u_0(i)$ for all $i\leq k$, we discuss
 it with two cases. Case 1: If $u (i) = {u_0} (i) =
 1$ for all $i \leq k$, we have $0 < \eta(u,1), \eta(u',1) < \frac
 {1}{k}$, thus $|\eta (u,1) - \eta (u_0,1)| < \frac {1}{k}$. Case 2:
 If there exists some state $i \leq k$ such that $u (i) = {u_0} (i) =
 0$, we denote the smallest such state as $i^*$ and we have $\eta (u,1) = \eta (u_0,1) =
 \frac{1}{i^*}$, thus $|\eta (u,1) - \eta (u_0,1)| = 0$. In
 summary, for any $\ep>0$, take $\hat k>1$ such that $\frac {1}{\hat k}<\ep$,
 thus $|\eta (u,1) - \eta (u_0,1)| <\ep$ for all $u\in O_{r^{\hat k}}(u_0)$. Therefore, $\eta (u,1)$ is continuous at
 $u_0$.


Finally, by Theorem \ref{thmopex}, the optimal stationary policy
exists. Actually, it is easy to verify that the optimal policy is
$u^* = (1,1,\dots,1,\dots)$ and the corresponding optimal cost is
$\eta^* = 0$.  \te
\end{example}

\begin{example}
(Example 8.10.2 in Puterman's book \citep{Puterman94})
 Consider an MDP with $\cS= \{1,2, \dots\}$. At state $i \in \cS$,
 there are two actions $1$ and $0$. If action $1$ is taken,
 then the state transits from $i$ to $i+1$ with probability $1$ and the reward is $f(i,1)=0$;
 if action $0$ is taken, then the state stays at
 $i$ with probability $1$ and the reward is $f(i,0)= 1 - \frac {1}{i}$. The Markov
 chain (under any policy) is denoted as $X(t)$, $t=0,1, \dots$. A stationary policy
 is denoted as a mapping $u: \cS \rightarrow \{0,1\} $.

 The performance measure for policy $u = (u(1), u(2),
 \dots ) $ with initial state $i$ is the long-run average reward
 as follows.
  \be    \label{defperf2}
  \eta (u,i) = \lim_{T \to \infty} \frac {1}{T} \mathbb E \left \{
  \sum_{t=0}^{T-1}
  f(X(t),u(X(t))) \Big | X(0) = i \right \} .
  \ee
 We set the initial state always as $X(0)=1$ and we wish to find an policy $u^*$ such that
 \[
 \eta (u^*,1) = \sup_{u \in \cal U} \{ \eta (u,1) \} .
 \]

 The discussion is the same as Example \ref{exmp}, except that
 $\eta (u,1)$ is NOT continuous at $u_0 = (1,1,  \dots, 1, \dots )$
 with $\eta (u_0,1) = 0$, while $\eta (u,1) \geq  1-\frac{1}{k}$ for any neighboring policy $u$ with $d(u,u_0)<r^k$.
 Therefore, an optimal stationary policy may not exist
 for this example. Actually, it is easy to verify that the optimal reward of this problem is $\eta^* = 1$.
 A history-dependent policy $u^*$ which uses action 0 $i$ times in state $i$, and then uses action 1 once,
 will yield a reward stream of $(0,0,\frac{1}{2},\frac{1}{2},0,\frac{2}{3},\frac{2}{3},\frac{2}{3},0,\frac{3}{4},\frac{3}{4},\frac{3}{4},\frac{3}{4},\dots)$.
 Thus, the history-dependent policy $u^*$ can reach the optimal reward $\eta^*=1$.
 However, any stationary deterministic policy yields
 possible rewards as either 0 or $1-\frac{1}{i}$, which cannot reach the
 optimal reward $\eta^*=1$.
 \te
\end{example}

\section{The $c/\mu$-Rule in Queueing Systems}\label{section_queue}
In this section, we show that, with the metric space defined by
(\ref{defdis}), the basic idea presented in Section
\ref{section_result} can be applied to a class of optimal scheduling
problems in queueing systems, called the $c/\mu$-rule problem, to
establish the existence of an optimal stationary policy.

\begin{figure}[htbp]
\centering
\includegraphics[width=0.6\columnwidth]{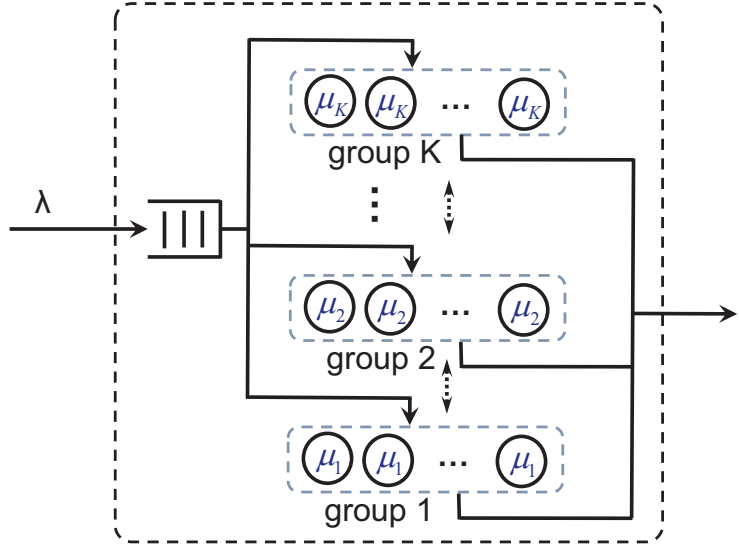}
\caption{The illustration of the on/off control of group-server
queues.}\label{fig_groupQ}
\end{figure}

The problem is about the on/off scheduling control of parallel
servers in a group-server queue. More details of the problem setting
can be referred to \citep{Xia2018} and we give a brief introduction
as follows. Consider a group-server queue with a single
infinite-size buffer and $K$ groups of parallel servers, as
illustrated by Fig.~\ref{fig_groupQ}. Customers are homogeneous and
customer arrival is assumed as a Poisson process with rate
$\lambda$. Arriving customers will go to the idle servers at status
`on'. If all the servers at status `on' are busy, the arriving
customer will wait in the buffer. Servers are providing service in
parallel and categorized into $K$ groups. Servers in the same group
are homogeneous in service rates and cost rates, while those in
different groups are heterogeneous. Group $k$ has $M_k$ servers with
service rate $\mu_k$ and cost rate $c_k$ per unit of time,
$k=1,2,\cdots,K$. Without loss of generality, we assume $\mu_1 \geq
\mu_2 \geq \dots \geq \mu_K$. The system cost includes two parts,
the operating cost of servers and the holding cost of customers. The
system state $n$ is the number of customers in the system. The state
space is denoted as $\cS = \{0,1,2,\dots\}$, which is countably
infinite. We can turn on or off servers dynamically to reduce the
system average cost. The action is the number of working servers at
each group, which is denoted as $a=(a_1,a_2,\cdots,a_K)$, where
$a_k$ is the number of working servers in group $k$ and $a_k \in
\{0,1, \cdots, M_k\}$. For any state $n \geq 1$, action space
$\mathcal A(n)$ is a subset of $\{1,\ldots,M_1\}\times
\{0,\ldots,M_2\}\times\cdots\times \{0,\ldots,M_K\}$, where $a_1
\geq 1$ is reasonable to guarantee the system ergodic. Define a
stationary policy as $u := (u(0),u(1),u(2),\dots)$, where
$u(n):=(u(n,1),u(n,2),\dots,u(n,K)) \in \mathcal A(n)$ is the action
at state $n$ and $u(n,k)$ is the number of working servers in group
$k$ at state $n$. The cost function at state $n$ under policy $u$ is
\begin{equation} \label{rewardf}
f(n,u) := h(n) + \sum_{k=1}^{K}c_k u(n,k),
\end{equation}
where $h(n)$ is the holding cost rate at state $n$. The system
long-run average cost under policy $u$ is defined as
\begin{equation}  \label{etau}
\eta (u) = \lim\limits_{T \rightarrow \infty} \frac{1}{T} \mathbb E
\left\{ \int_{t=0}^{T} f(n(t),u)dt \right\},
\end{equation}
where $n(t)$ is the system state at time $t$. The optimal average
cost is $\eta^* = \inf\limits_{u} [\eta(u)]$. We aim at finding the
optimal stationary policy $u^*$ which achieves the optimal average
cost, i.e., $\eta(u^*) = \eta^*$, where $u^* \in \mathcal U$ and
$\mathcal U$ is the stationary policy space. In \citep{Xia2018}, it
is shown that the optimal policies (if one exists) follow the so
called $c/\mu$-rule: Servers in the group with smaller values of
$c/\mu$ should be turned on with higher priority. Here, we want to
verify that an optimal stationary policy does exist for this problem
with countable states.

It is natural to assume that the holding cost $h(n)$ is increasing
in $n$; and thus, under optimal policies the queue should be
\emph{ergodic}. So we assume that $\{n(t)\}$  is  ergodic (under
each policy in $\cU$) with a unique steady-state distribution $\pi
(n,u)$, $n=0,1, \cdots$, $u \in \cU$, and the long-run average
(\ref{etau}) does not depend on the initial state.

Since our queue is a birth-death process, we
can derive the steady-state distribution as below.
\begin{eqnarray}\label{eq_pi}
\pi(n,u) = \frac{1}{1+G(u)}\prod_{l=1}^{n}\frac{\lambda}{u(l)\mu},
\qquad n \geq 1,
\end{eqnarray}
where $\mu = (\mu_1, \cdots, \mu_K )^T$, and
\begin{equation}
u(l)\mu := \sum_{k=1}^{K} u(l,k) \mu_k,
\end{equation}
and
\begin{equation} \label{gun}
G(u) := \sum_{n=1}^{\infty} \prod_{l=1}^{n}\frac{\lambda}{u(l)\mu}.
\end{equation}
The queue is stable if and only if $G(u) < \infty$ which also
indicates
\begin{equation}\label{eq7}
\lim\limits_{n \rightarrow \infty} \sum_{m=n}^{\infty}
\prod_{l=1}^{m}\frac{\lambda}{u(l)\mu} =0.
\end{equation}
The ergodicity of the system under a policy $u$ can indicate a
necessary condition: $u(n)\mu \neq 0$ for all $n \geq 1$. The
stability of the system can be guaranteed by a sufficient condition:
there exists an $\bar n$ such that $u(n)\mu
> \lambda$ for all $n > \bar n$.

For an ergodic policy $u \in \cU$, under suitable condition, the long-run average (\ref{etau}) equals
\be  \label{etast}
\eta (u) = \sum_{n=0}^\infty \pi (n,u)f(n,u) .
\ee

For the analysis here, we need to make the following assumption:
\begin{assumption} \label{ass0}
The normalizing factor $G(u)$ in (\ref{gun}) (equivalently, the limit in (\ref{eq7})) and the performance limit (\ref{etast}) converge
uniformly in $\cU$.
\end{assumption}

%
%

We use the metric definition (\ref{defdis}) to quantify the distance
between any two policies $u_1$ and $u_2$. In what follows, we will
prove that when the two policies $u$ and $u'$ are infinitely close,
their performance measures $\eta (u)$ and $\eta (u')$ are also
infinitely close to each other. Denote the two policies by $u =
(u(0), u(1), \cdots, u(n), \cdots )$ and $u' = (u'(0), u'(1),
\cdots, u'(n), \cdots )$.
By Lemma \ref{lemduout}, we assume that
\be  \label{uln}
u(l)=u'(l),~~for ~l=0,1, \cdots, n.
\ee
which means that $u'\in O_{r^n}(u)$.

First, we compare the difference of the normalization factors
$1+G(u)$ and $1+G(u')$ of these two policies. We have
\begin{eqnarray} \label{onepd}
\frac{1+ G(u)}{1+G(u')} &=& \frac{1+\sum_{m=1}^{\infty}
\prod_{l=1}^{m}\frac{\lambda}{u(l)\mu}}{1+\sum_{m=1}^{\infty}
\prod_{l=1}^{m}\frac{\lambda}{u'(l)\mu}} \nonumber\\
&=& \frac{  \Big ( 1+\sum_{m=1}^{n}
\prod_{l=1}^{m}\frac{\lambda}{u(l)\mu} \Big ) +\sum_{m=n+1}^{\infty}
\prod_{l=1}^{m}\frac{\lambda}{u(l)\mu}}{ \Big ( 1+\sum_{m=1}^{n}
\prod_{l=1}^{m}\frac{\lambda}{u(l)\mu} \Big ) +\sum_{m=n+1}^{\infty}
\prod_{l=1}^{m}\frac{\lambda}{u'(l)\mu}} \nonumber\\
&=& \frac{1+\frac{\sum_{m=n+1}^{\infty}
\prod_{l=1}^{m}\frac{\lambda}{u(l)\mu}}{1+\sum_{m=1}^{n}
\prod_{l=1}^{m}\frac{\lambda}{u(l)\mu}}}
{1+\frac{\sum_{m=n+1}^{\infty}
\prod_{l=1}^{m}\frac{\lambda}{u'(l)\mu}}{1+\sum_{m=1}^{n}
\prod_{l=1}^{m}\frac{\lambda}{u(l)\mu}}} <
1+\frac{\sum_{m=n+1}^{\infty}
\prod_{l=1}^{m}\frac{\lambda}{u(l)\mu}}{1+\sum_{m=1}^{n}
\prod_{l=1}^{m}\frac{\lambda}{u(l)\mu}}
\nonumber\\
&<& 1 + \sum_{m=n+1}^{\infty} \prod_{l=1}^{m}\frac{\lambda}{u(l)\mu}
= 1 + \delta(n,u),
\end{eqnarray}
where
\begin{equation}
\delta(n,u) := \sum_{m=n+1}^{\infty}
\prod_{l=1}^{m}\frac{\lambda}{u(l)\mu}.
\end{equation}
Similarly, we can also have
\begin{eqnarray}   \label{onemdp}
\frac{1+ G(u)}{1+G(u')} &=& \frac{1+\sum_{m=1}^{\infty}
\prod_{l=1}^{m}\frac{\lambda}{u(l)\mu}}{1+\sum_{m=1}^{\infty}
\prod_{l=1}^{m}\frac{\lambda}{u'(l)\mu}} =
\frac{1+\frac{\sum_{m=n+1}^{\infty}
\prod_{l=1}^{m}\frac{\lambda}{u(l)\mu}}{1+\sum_{m=1}^{n}
\prod_{l=1}^{m}\frac{\lambda}{u(l)\mu}}}
{1+\frac{\sum_{m=n+1}^{\infty}
\prod_{l=1}^{m}\frac{\lambda}{u'(l)\mu}}{1+\sum_{m=1}^{n}
\prod_{l=1}^{m}\frac{\lambda}{u(l)\mu}}} \nonumber\\
&>& \frac{1} {1+\frac{\sum_{m=n+1}^{\infty}
\prod_{l=1}^{m}\frac{\lambda}{u'(l)\mu}}{1+\sum_{m=1}^{n}
\prod_{l=1}^{m}\frac{\lambda}{u(l)\mu}}}
> \frac{1} {1+ \sum_{m=n+1}^{\infty}
\prod_{l=1}^{m}\frac{\lambda}{u'(l)\mu} }\nonumber\\
&>& 1 - \sum_{m=n+1}^{\infty}
\prod_{l=1}^{m}\frac{\lambda}{u'(l)\mu} = 1-\delta(n,u'),
\end{eqnarray}
where
\begin{equation}
\delta(n,u' ) := \sum_{m=n+1}^{\infty}
\prod_{l=1}^{m}\frac{\lambda}{u'(l)\mu}.
\end{equation}
%
%
Therefore, we have
\begin{eqnarray}
1-\delta(n,u' ) < \frac{1+ G(u)}{1+G(u')} < 1+\delta(n,u).
\end{eqnarray}
Let $\sigma(n,u,u')$ be determined by
\begin{equation}\label{eq_G2}
\frac{1+ G(u)}{1+G(u')} = 1+\sigma(n,u,u').
\end{equation}
Then,
\begin{equation}\label{eq_sigma}
-\delta(n,u')<\sigma(n,u,u')<\delta(n,u) .
\end{equation}
%
With (\ref{eq_pi}), (\ref{uln}), and (\ref{eq_G2}), the steady-state distributions
under these two policies $u$ and $u'$ have the following relation.
\begin{equation}\label{eq13}
\pi(m,u') = (1+\sigma(n,u,u'))\pi(m,u),  \qquad m=0,1,\dots,n.
\end{equation}

Next, we study the difference between the associated long-run
average costs $\eta$ under policies $u$ and $u'$. The cost functions
are denoted by $f(m,u)$ and $f(m,u')$, respectively. By
(\ref{rewardf}) and (\ref{uln}), we have $f(m,u) = f(m,u')$ for
$m=0,1,\dots,n$. Therefore, we have
\begin{eqnarray}
&&\eta(u' ) - \eta(u) \nonumber\\
&=& \sum_{m=0}^{\infty}[\pi(m,u')f(m,u' ) - \pi(m,u)f(m,u)] \nonumber\\
&=& \sum_{m=0}^{n}[\pi(m,u')f(m,u' ) - \pi(m,u)f(m,u)] + \sum_{m=n+1}^{\infty}[\pi(m,u')f(m,u') - \pi(m,u)f(m,u)] \nonumber\\
&=& \sum_{m=0}^{n}[\pi(m,u') - \pi(m,u)]f(m,u) +
\sum_{m=n+1}^{\infty}[\pi(m,u')f(m,u') - \pi(m,u)f(m,u)] .  \nonumber
\end{eqnarray}
Applying (\ref{eq13}), we have
\begin{equation}\label{diffeta}
\eta(u') - \eta(u) = \sigma(n,u,u')\sum_{m=0}^{n}\pi(m,u)f(m,u) +
\sum_{m=n+1}^{\infty}[\pi(m,u')f(m,u') - \pi(m,u)f(m,u)].
\end{equation}
%

Now we are ready to prove the continuity of $\eta (u)$ in the metric
space $\cU$ with metric (\ref{defdis}). With (\ref{eq7}), we have
\begin{equation*}
\lim\limits_{n \rightarrow \infty} \delta(n,u) = 0, \qquad
\lim\limits_{n \rightarrow \infty} \delta(n,u') = 0.
\end{equation*}
Let $\ep>0$ be any small number. Under Assumption~\ref{ass0}, by the
uniformity of $G(u)$ in (\ref{gun}) and (\ref{eq7}), there exists a
large integer $N_1$ such that if $n>N_1$, we have $\delta (n,u) <
\ep$ for any $u \in \cU$. By (\ref{eq_sigma}), we have
\[
|\sigma (n,u,u')|< \ep, \quad \forall u,u' \in \cU.
\]
Next, because (\ref{etast}) converges, there is a large integer $N_2
$ such that
\begin{equation*}
\Big|\sum_{m=0}^{n}\pi(m,u)f(m,u) \Big|< |\eta(u)| + 1, \quad
\forall n>N_2.
\end{equation*}
Furthermore, under Assumption \ref{ass0}, by the uniformity of the
convergence of (\ref{etast}), there is a large integer $N_3$ such
that
\begin{equation*}
\Big|\sum_{m=n+1}^{\infty}[\pi(m,u')f(m,u') - \pi(m,u)f(m,u)] \Big|
< 2 \ep, \quad \forall n>N_3 \ {\rm and} \ u,u'\in  \cU.
\end{equation*}
Finally, let $N^*:=\max \{N_1, N_2, N_3 \}$. Then, by (\ref{diffeta}) and Lemma~\ref{lemduout},  we
have
\begin{eqnarray}\label{diffop1}
|\eta (u) - \eta (u')| &\leq&
|\sigma(n,u,u')|\Big|\sum_{m=0}^{n}\pi(m,u)f(m,u)\Big| +
\Big|\sum_{m=n+1}^{\infty}[\pi(m,u')f(m,u') - \pi(m,u)f(m,u)] \Big|
\nonumber\\
&<& [|\eta(u)|+ 3] \ep, \ \ \ \ {\rm for \ all \ }       u'\in
O_{r^{N^*}}(u).
\end{eqnarray}
Since $\eta(u)$ is bounded, we conclude that $\eta (u)$ is
continuous at $u$ in the metric space. Therefore, the existence of
optimal stationary policy $u^*$ for this $c/\mu$-rule problem
directly follows by Theorem~\ref{thmopex}. \te

\noindent\textbf{Remark~2.} The condition of uniform convergence in
Assumption~\ref{ass0} is easy to validate in queueing systems. For
example, we can set the condition for the control of our
group-server queues as follows: \textcircled{\oldstylenums{1}} there
exists a constant $\tilde{n}$ such that for any $n
> \tilde{n}$, every feasible action $u(n) \in \mathcal A(n)$ always satisfies $u(n)
\mu > \lambda$. Therefore, we define $\rho_0 := \max_{u(n)\in
\mathcal A(n),n>\tilde{n}}\{ \frac{\lambda}{u(n)\mu} \} <1$. We
directly have $G(u)\leq
\sum_{n=1}^{\tilde{n}}\prod_{l=1}^{n}\frac{\lambda}{u(l)\mu}+\sum_{n=\tilde{n}+1}^{\infty}\rho_0^n
< \infty$, which indicates that the queueing system is stable and
the normalizing factor $G(u)$ in (\ref{gun}) converges uniformly in
$u\in \cU$. Compared with (\ref{eq_pi}), we further define a pseudo
probability $\tilde{\pi}(n,u) := \frac{1}{1+G(u)}\rho_0^n$.
Obviously, we always have $\tilde{\pi}(n,u) \geq \pi(n,u)$ for any
policy $u$ and $n > \tilde{n}$. Thus, for the performance limit
(\ref{etast}), we have $|\eta (u)| \leq \sum_{n=0}^{\tilde n}
\pi(n,u)|f(n,u)| + \sum_{n=\tilde{n}+1}^{\infty}
\tilde{\pi}(n,u)|f(n,u)| = \sum_{n=0}^{\tilde n} \pi(n,u)|f(n,u)| +
\frac{1}{1+G(u)} \sum_{n=\tilde{n}+1}^{\infty} \rho_0^n |f(n,u)|$,
where the first part is always finite and we only need to guarantee
the second part bounded. Thus, \textcircled{\oldstylenums{2}} any
cost function $|f(n,u)|$ polynomially increasing to infinity along
with $n$ will be controlled by the exponential factor $\rho_0^n$.
Therefore, with \textcircled{\oldstylenums{1}} and
\textcircled{\oldstylenums{2}}, we can easily validate
Assumption~\ref{ass0} that $G(u)$ and $\eta(u)$ converge uniformly,
and thus an optimal stationary policy exists. More specifically, for
the cost function (\ref{rewardf}), we have
$f(n,u)=h(n)+\sum_{k=1}^{K}c_k u(n,k)$, where the operating cost
$\sum_{k=1}^{K}c_k u(n,k)$ is obviously bounded and the holding cost
$h(n)$ can be unbounded. From the above analysis, we can see that
$f(n,u)$ can be unbounded from both below and above sides. For
example, we can set $h(n)=(-1)^n \cdot n$, which is unbounded both
below and above while satisfies our condition
\textcircled{\oldstylenums{2}}. However, this kind of cost function
may not be handled by other methods in the literature
\citep{Borkar1994} because the cost function thereof is required to
be unbounded from  below. This is also one of the advantages of our
method in this paper.

We have demonstrated the applicability of Theorem~\ref{thmopex} for
proving the existence of optimal stationary policies in a scheduling
problem of queueing systems. In the next section, we further show
that this approach also applies to more general cases.

\section{More General Cases}\label{section_extension}
In general, we consider a continuous-time MDP with a countable state
space denoted as $\cS = \{ 0,1,   \dots \}$. Let $\pi (i,u)$ be the
steady-state probability of state $i\in \cS$ under given policy
$u\in \mathcal U$,  and $q^a(i,j)$ be the transition rate from state
$i$ to $j$ under action $a\in \mathcal A(i)$, $i,j \in \cS$.
Obviously, we have $q^a(i,j)\geq 0$ for $i \neq j$ and
$q^a(i,i)=-\sum_{j\in\mathcal S, j \neq i} q^a(i,j) \leq 0$, where
$|q^a(i,i)|$ can be understood as the total rates transiting out
from state $i$ if action $a$ is adopted. Then we know that the
steady-state probabilities $\pi(i,u)$'s must satisfy the following
equations.
\begin{align}
& \sum_{j=0}^\infty \pi (j,u) q^{u(j)}(j,i) = 0, \quad i \in \cS, \label{pii} \\
& \sum_{i=0}^\infty \pi (i,u) = 1 ,    \label{sumpi}
\end{align}
where (\ref{sumpi}) is called a normalization equation. Given a
policy $u\in \mathcal U$, any sequence $\nu(i,u)\geq 0$ (depending
on $u$), $i \in \cS$, that satisfies
\begin{equation} \label{eqnu}
\sum_{j=0}^\infty \nu (j,u) q^{u(j)}(j,i) = 0, \ \forall i \in \cS \
\mbox{ and} \quad \sum_{i=0}^\infty \nu (i,u) < \infty ,
\end{equation}
is called an {\it un-normalized steady-state vector}. From
(\ref{eqnu}), we have
\[
\pi (i,u) = \frac {\nu (i,u) } {\sum_{i=0}^\infty \nu (i,u) }, ~~i \in
\mathcal S ,
\]
is the steady-state probability.

In the rest of the paper, it is more convenient to deal with the
un-normalized vector because it does not contain the denominator.
Moreover, it is convenient to set $\nu (0,u)=1$ to obtain an
un-normalized probability.

First, we make the following assumptions to simplify the problem
setting.
\begin{assumption}  \label{ass1}
\begin{enumerate}
\item [(a)] $q^a(i,j)$ is bounded, i.e., $|q^a(i,j)| < \Lambda $, for all $i,j\in \mathcal S, a\in \mathcal A(i)$.
\item [(b)] There is an integer $M>0$ such that $q^a(j,i)=0$, for
all $j >i+M$, $i \in \mathcal S$ and $a\in \mathcal A(i)$.
\end{enumerate}
\end{assumption}
Assumption~\ref{ass1}(a) indicates that the transition rate from any
state $i$ has an upper bound $\Lambda$, which is reasonable for most
cases in practice. Assumption~\ref{ass1}(b) means that the
transition rate from state $j$ back to $i$ is 0 if state $j$ is far
away from state $i$. This assumption is also reasonable in many
practical systems, especially it is usually true for queueing
systems since state $j$ always transits back only to state $j-1$
caused by a service completion event.

Given any $u\in \cU$, at a state $i \in \cS$, we may take an action
denoted by $u(i)$, which determines the value of $q^{u(i)}(i,j)$, $j
\in \cS$. Then $u:= (u(0), u(1), \cdots )$ denotes a policy. Let
$\cU$ be the space of all policies. The steady-state probability at
state $i$ is denoted by $\pi(i,u)$, which depends on policy $u$. The
reward or cost function at state $i$ with action $u(i)$ is denoted
by $f(i,u(i))$. We assume that the Markov processes under all
policies in $\cU$ are ergodic and the long-run average performance
under policy $u$ is \be \label{etaepf} \eta (u): = \sum_{i=0}^\infty
\pi (i,u) f(i,u(i)). \ee

Denoting $\nu(i,u)$ as the un-normalized steady-state vector
satisfying (\ref{eqnu}) under policy $u$, we give one more
assumption as follows (cf. Assumption \ref{ass0}).
 \begin{assumption} \label{ass2}
  $\sum_{i=0}^N \nu(i,u) $, with $\nu(0,u)=1$, converges uniformly in
  $\cU$ as $N \to \infty$, and $\sum_{i=0}^N \pi(i,u) f(i,u(i))$
  converges uniformly in $\cU$, as $N \to \infty$.
 \end{assumption}

Assumption~\ref{ass2} holds for many Markov systems, especially when
the system is stable under the neighborhood of policies. In fact, it
holds if there is a sequence, denoted as $\overline \nu (i)$,
$i=0,1, \dots$, such that $\nu(i,u) \leq \overline \nu (i)$ and
$\sum_{i=0}^\infty \overline\nu (i)< \infty$.

\begin{example}
Consider a controlled $M/M/1$ queue with arrival rate $\lambda (i,u)$ and service
rate $\mu (i,u)$ (under a given control policy $u$)  when the number of customers is $i$, $i\in \cS = \{0,
1, \dots, \}$. Let $X(t) \in \cS$ be the Markov process of the
queue. The un-normalized steady-state vector is $\nu (i,u) =
\prod_{l=0}^i \frac {\lambda(l,u)}{\mu(l,u)}$.
 The process is stable if
 \[
 \sum_{i=0}^\infty \nu (i,u) = \sum_{i=0}^\infty
 \prod_{l=0}^i \frac {\lambda(l,u)}{\mu(l,u)} < \infty .
 \]
Therefore, Assumption \ref{ass2} is the same as Assumption \ref{ass0},
and if there is a bound $\overline \gamma <1$ and state $i^*$ such that
$\frac {\lambda(i,u)}{\mu(i,u)} < \overline \gamma $ for all
policies $u$ and states $i\geq i^*$, then Assumption \ref{ass2} holds. \te
\end{example}

Now, let us understand the role of Assumptions~\ref{ass1} and
\ref{ass2}. For any integer $N>0$, we consider the first $K$
equations in (\ref{eqnu}), where $K>N$. Given any $u\in \cU$, by
Assumption~\ref{ass1}(b), the summation in (\ref{eqnu}) is over only
finitely many states, resulting in
\begin{equation}\label{eqshort}
\sum_{j=0}^{i+M} \nu (j,u)q^{u(j)}(j,i) = 0,
 ~~~i=0,1, \dots, K,
\end{equation}
which can be further rewritten as
\begin{equation}\label{twoterm0}
\sum_{j=0}^{K} \nu (j,u) q^{u(j)}(j,i)
  + \sum_{j=K+1}^{i+M} \nu (j,u) q^{u(j)}(j,i) = 0,
 ~~~i=0,1, \dots, K.
\end{equation}
 For   (\ref{twoterm0}), the last summation is nonzero only if
$i+M >K$. Thus, only the last $M$ equations in (\ref{twoterm0})
contain nonzero terms of the last summation, whose values are small
enough to be ignored, as shown by the following analysis.

 For any $\ep >0$ and
$N>0$, by Assumption  \ref{ass2},  there is a large enough $K$ such that
\begin{equation}\label{sumnuk}
\sum_{i=K+1}^\infty \nu (i,u) < \frac {\ep}{N}, \ \ \ {\rm for \ all}  \  u\in \cU.
\end{equation}
By Assumption~\ref{ass1}, the last summation of (\ref{twoterm0}) can
be written as
\begin{equation}\label{sumnuk2}
\sum_{j=K+1}^{i+M} \nu (j,u) q^{u(j)}(j,i) < \sum_{j=K+1}^{\infty}
\nu (j,u) q^{u(j)}(j,i) < \frac {\ep}{N} \Lambda = O(\frac
{\ep}{N}), \ \ \ {\rm for \ all}  \  u\in \cU.
\end{equation}
Substituting the above result into (\ref{twoterm0}), we  see that solving (\ref{twoterm0}) becomes solving the following
equations
\begin{eqnarray}\label{twoterm3}
0 &=& \sum_{j=0}^{K} \nu (j,u) q^{u(j)}(j,i),
 ~~~~~~~~~~~~~~~ i=0,1, \dots, K-M, \nonumber \\
0 &=& \sum_{j=0}^{K} \nu (j,u) q^{u(j)}(j,i) + O(\frac {\ep}{N}),
 ~~~i=K-M+1, \dots, K,
\end{eqnarray}
where we have $K+1$ variables and $K+1$ linear equations. Thus, the
variables $\nu(i,u)$'s can be solved and we state the results as
(\ref{muinep}) in the following lemma, where $F_i (q^{u(j)}(j,k); \
j,k=0,1, \dots, K )$ denotes a function $F_i(\cdot)$ with variables
$q^{u(j)}(j,k)$, $i=0,1, \dots, K$.


\begin{lemma}  \label{lemnu}
Under Assumptions ~\ref{ass1} and \ref{ass2}, for any policy
$u\in\cU$, integer $N>0$, and small number $\ep >0$, there exists an
integer $K>0$ such that
\begin{equation} \label{muinep}
\nu (i,u) = F_i (q^{u(j)}(j,k); \   j,k=0,1, \dots, K ) + \kappa_i
(N) ,~~~~i=0,1, \dots, N,\footnote{In fact, this equation holds for
$i=0,1, \dots, K$, $K>N$, but to prove Theorem \ref{thm2}, we only
need it for the first $N$ $\nu (i,u)$'s.}
\end{equation}
and $\kappa_i (N) < \frac {\ep}{N}$. In words, we say that roughly
for any finite $N$, $\nu (0,u) , \dots, \nu (N,u)$ depend only on
the transition rates among finitely many states. The functions
$F_i$, $i=0,1,\cdots, N$, are the same for any policy $u' \in
O_{r^K}(u)$.
\end{lemma}

Note that we can set
$\nu(0,u)=1$ for solving (\ref{twoterm3}) since $c \nu$ is also a
solution to (\ref{twoterm3}) for any feasible solution $\nu$, where
$c$ is a constant. Moreover, ignoring the term of $O(\frac
{\ep}{N})$, (\ref{twoterm3}) is a set of linear equations determined
by the values of $\{q^{u(j)}(j,k); \ j,k=0,1, \dots, K \}$.
Therefore, for any two  policies $u'$ and $u$ such that $u'(i)=u(i)$
for all $0\leq i\leq K$, $F_i$'s take the same form for such
policies, $i=0,1, \cdots, K$.

With Assumptions~\ref{ass1} and \ref{ass2}, we can further extend
the existence condition of optimal stationary policies in
Theorem~\ref{thmopex} and derive the following theorem.
\begin{theorem}\label{thm2}
Under Assumptions~\ref{ass1} and \ref{ass2}, there exists an optimal
stationary policy for the average cost MDP with a countable state
space.
\end{theorem}
{\it Proof:} 
 %
 %
 Let $N>0$ be any integer and $\ep>0$ be any small number.
 Consider any two policies $u$ and $u'$, which determine the corresponding transition rates $q(i,j):=q^{u(i)}(i,j)$ and $q'(i,j):=q^{u'(i)}(i,j)$, as well as the steady-state vectors $\nu(i)$ and $\nu'(i)$, respectively.
 By Lemma~\ref{lemnu} and Assumption~\ref{ass2}, if $K$ is large enough,
 then we have
 \[
\nu ' (i) = F_i (q'(j,k);\ j,k=0,1, \dots, K ) + \kappa'_i (N)
,~~~~i=0,1, \dots, N,
 \]
 and
  \[
\nu (i) = F_i (q(j,k);\ j,k=0,1, \dots, K ) + \kappa_i (N)
,~~~~i=0,1, \dots, N,
 \]
where $ \kappa '_i (N) < \frac {\ep}{N}$ and $ \kappa _i (N) < \frac
{\ep}{N}$.

By Lemma~\ref{lemduout}, if $u$ and $u'$ are close enough such that
$d(u, u' )<r^K$, then $u (i) = u' (i)$ for all $i<K$. This means
$q(i,j) = q'(i,j)$ for all $i<K$ and $j=0,1, \dots$. Therefore, we
have
\begin{equation} \label{nudiff}
\nu ' (i) = \nu (i) + \kappa
'_i (N) - \kappa_i (N) , ~~i=0,1, \dots, N.
\end{equation}

The rest analysis is similar to (\ref{onepd})--(\ref{eq13}). First,
we have
\[
\pi (i) = \frac {\nu (i)}{\sum_{j=0}^\infty \nu (j) } ,
\]
and
\begin{eqnarray}  \label{pippi}
\pi' (i) &=& \frac {\nu ' (i)}{\sum_{j=0}^\infty \nu ' (j) } \nonumber \\
&=&  \frac {\sum_{j=0}^\infty \nu (j) }{\sum_{j=0}^\infty \nu ' (j)
} \Big \{ \frac {\nu (i) + \kappa_i (N) - \kappa'_i (N) }
{\sum_{j=0}^\infty \nu  (j) } \Big \}, \quad i=0,1,\dots,N.
\end{eqnarray}
With (\ref{nudiff}), we have
\begin{eqnarray}
\frac {\sum_{j=0}^\infty \nu (j) }{\sum_{j=0}^\infty \nu ' (j) } &=&
\frac {\sum_{j=0}^N \nu (j) + \sum_{j=N+1}^\infty \nu (j) }
{\sum_{j=0}^N \nu ' (j) + \sum_{j=N+1}^\infty \nu ' (j) } \nonumber \\
&=& \frac {\sum_{j=0}^N \nu (j) + \sum_{j=N+1}^\infty \nu (j) }
{\sum_{j=0}^N \nu  (j) + \sum_{j=N+1}^\infty \nu ' (j) +
\sum_{j=0}^N [\kappa'_i (N) - \kappa_i (N) ]} . \nonumber
 \end{eqnarray}
If $K$ is large enough (i.e., $d(u,u')$ is small enough), it holds
 \[
 \Big | \sum_{j=0}^N [\kappa'_i (N) - \kappa_i (N) ] \Big | < 2 \ep .
 \]
Therefore,
\begin{align}
\frac {\sum_{j=0}^\infty \nu (j) }{\sum_{j=0}^\infty \nu ' (j) } & =
\frac {1 + \frac {\sum_{j=N+1}^\infty \nu (j)} {\sum_{j=0}^N \nu (j)
} } {1 + \frac {\sum_{j=N+1}^\infty \nu ' (j)} {\sum_{j=0}^N \nu (j)
} + \frac {\sum_{j=0}^N [\kappa'_i (N) - \kappa_i (N) ]}
{\sum_{j=0}^N \nu(j)}} \nonumber \\
& =
\frac {1 + \frac {\sum_{j=N+1}^\infty \nu (j)} {\sum_{j=0}^N \nu (j)
} } {1 + \frac {\sum_{j=N+1}^\infty \nu ' (j)} {\sum_{j=0}^N \nu (j)
} + \ep (N, u,u') } , \nonumber
\end{align}
with $|\ep (N,u,u')| := \Big|\frac {\sum_{j=0}^N [\kappa'_i(N) -
\kappa_i(N) ]} {\sum_{j=0}^N \nu(j)}\Big| < \Big|\frac {\sum_{j=0}^N
[\kappa'_i(N) - \kappa_i(N) ]} {1} \Big|  < 2\ep$, where we use the
preset condition $\nu(0)=1$.

The rest proof follows the same procedure as
(\ref{onemdp})--(\ref{eq13}). First, as in (\ref{eq13}), we can
derive
\begin{equation}
\pi'(i) = (1+\sigma(N,u,u'))\pi (i), \quad i=1,2,\dots,N,
\end{equation}
where $|\sigma (N,u,u')|< 3 \ep$ (with a large $N$ such that
$\sum_{i=N+1}^\infty \nu (i) < \ep$), when $K$ is large enough.
Then, similar to (\ref{diffeta}), we have
\begin{eqnarray}  \label{diffeta2}
&&\eta(u') - \eta(u) \nonumber\\
&=&  \sigma(N,u,u')\sum_{m=0}^{N}\pi(i,u)f(i,u(i) ) +
\sum_{i=N+1}^{\infty}[\pi(i,u' )f(m,u'(i) ) - \pi(i,u)f(i,u (i) )] .  
\end{eqnarray}

Now we are ready to prove the continuity of $\eta (u)$ in the metric
space $\cU$ with metric (\ref{defdis}). Let $\ep>0$ be any small
number. First, as discussed above, under Assumptions \ref{ass1} and
\ref{ass2}, by the uniformity of $\sum_{i=0}^\infty \nu(i)$, there
is a large integer $N_1$ such that if $n>N_1$, we have $|\sigma
(N_1,u,u')|< 3\ep$ for any $u$ and $u'$. Next, because
$\sum_{i=0}^\infty \pi (i,u) f(i, u(i)) $ converges, there is an
$N_2 $ such that $|\sum_{i=0}^{N}\pi(i,u)f(i,u(i)) | < |\eta
(u)|+1$, for all $N>N_2$. Furthermore, under Assumption \ref{ass2},
by the uniformity of the convergence of (\ref{etaepf}), there is a
large $N_3$ such that for all $n>N_3$, it holds
\[
\Big|\sum_{i=n+1}^{\infty}[\pi(i,u')f(i,u') - \pi(i,u)f(i,u)] \Big|
< 2 \ep .
\]
Therefore, by (\ref{diffeta2}) and Lemma~\ref{lemduout}, for  $\hat
N:=\max \{N_1, N_2, N_3 \}$, we have
\begin{equation}\label{diffop}
|\eta (u) - \eta (u')| < [3|\eta(u)|+ 5] \ep, \ \ \ \forall \ u'\in
O_{r^{\hat N}}(u).
\end{equation}
Thus, $\eta (u)$ is continuous at $u$ in the metric space, and then
the existence of optimal stationary policy $u^*$ follows from
Theorem~\ref{thmopex}. \te
%

In summary, we have extended the existence condition of optimal
stationary policies for average MDPs with countable state space from
Theorem~\ref{thmopex} for the $c/\mu$-rule problem to
Theorem~\ref{thm2} for the more general case.
As stated by Assumptions~\ref{ass1} and \ref{ass2}, if the system
has bounded and limited-distance backward transition rates, and with
the uniformity of the convergence of the un-normalized probabilities
and the performance sequences, the existence of optimal stationary
policies can be guaranteed by Theorem~\ref{thm2}. The theorem may be
easily verified in practice, especially for queueing systems, as
demonstrated in the aforementioned examples.

%
%

\section{Conclusion}\label{section_conclusion}
In this paper, we derive the existence conditions of optimal
stationary policies for countable state MDPs with long-run average
criterion. By defining a suitable metric on the policy space forming
a compact metric space, the existence condition can be guaranteed by
proving the continuity of the long-run average cost as a function in
the policy space under the metric. With some assumptions
 on the transition rates and the uniformity of the convergence of the
 un-normalized probabilities of the processes, the existence of the
 optimal policies can be proved for the MDPs with countable states
 in a general form. Compared with other conditions studied in the
 literature, the condition in this paper may be easier to verify when
 applied to practical MDP problems, especially in queueing systems. Some
examples are studied to illustrate the applicability of our results.
Future research topics may include the extensions to MDPs with other
criteria, such as the discounted ones.

\section*{Acknowledgement}
The first author would like to thank Prof. Peter W. Glynn at
Stanford University for his comments, which partly initiate the work
of this paper.

This work was supported in part by the National Natural Science
Foundation of China (11931018, 61573206).

\end{document}